# More counterexamples to Lagrangian Poincaré recurrence in dimension four

Joel Schmitz

**Abstract**

In earlier work, we constructed counterexamples to Lagrangian Poincaré recurrence for many toric symplectic four manifolds. Here we provide a few more examples extending the family of counterexamples to include all non-monotone toric symplectic four manifolds.

*This work has been incorporated into [3, Section 3.3].*

This note is a follow up on [2]. The goal is again to find counter examples to Lagrangian Poincaré recurrence in dimension four, that is to construct a Hamiltonian diffeomorphism $\psi$ and a Lagrangian $L$ such that $\psi^n(L) \cap L = \emptyset$ for all $n \in \mathbb{N}$. See [3] for a short discussion of known results related to Lagrangian Poincaré recurrence.

As mentioned in [2, Remark 3.2], the construction of [2] can be carried out in every symplectic toric four manifold except:

1. The five monotone symplectic toric four manifolds.
2. $(S^2 \times S^2, \omega_{a,b})$ with $a \neq b$.
3. The one or two-fold blow-up of $(S^2 \times S^2, \omega_{a,b})$ of size $c = \min\{a,b\}/2$.

In this note we provide a slightly different construction that also provides counterexamples for cases 2 & 3, and we give a few more details on the construction hinted at in [2, Remark 3.2]. We will assume that the reader is familiar with the techniques developed in [3].

We first give the construction in the special case of non-monotone $S^2 \times S^2$ in Section 1, than repeat the construction in the general case in Section 2.

## 1 Non-monotone $S^2 \times S^2$

**The base space.** Let $M, w > 0$ and take $\pi_0 : (S^2 \times S^2, \omega_{2M+w,2M}) \to B_0$ to be the almost toric fibration with the nodal integral affine base space $B_0$ as in Figure 1 a). The symplectic form $\omega_{2M+w,2M}$ assigns areas $2M + w$ and $2M$ to the factors of $S^2 \times S^2$ respectively. Let $\mathscr{F} : B_0 \to [0, M]$ be the *height function* ([3, Section 3.1]) of $B_0$, that is the integral affine distance to $\partial B_0$. The level sets of $\mathscr{F}$ are drawn in grey in Figure 1 a). Let $\varepsilon > 0$. We assume that the four nodes in Figure 1 a) have height $\mathscr{F} > M - \varepsilon$. In particular this means that for $h \in [0, M - \varepsilon)$ the level sets $\mathscr{F}^{-1}(h)$ are straight lines and circles in $B_0$.

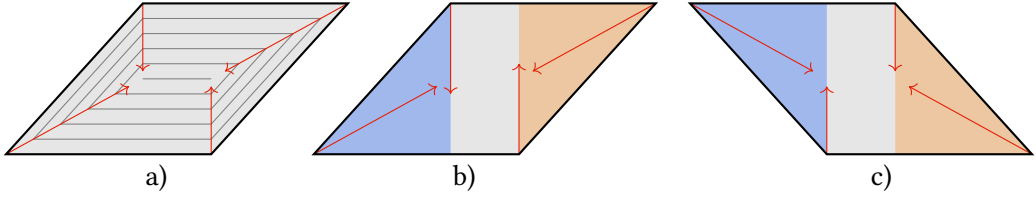

Figure 1: The construction for non-monotone $S^2 \times S^2$.

**The nodal tangle.** Modify $B_0$ by nodal slides to obtain the nodal integral affine surface $B_1$ pictured in Figure 1 b). Denote the resulting nodal tangle by $B \times [0, 1]$. We may choose a different nodal chart than Figure 1 b) for $B_1$: Cut off the blue and orange shaded triangles in Figure 1 b), apply the shear matrix $\begin{pmatrix} 1 & 0 \\ -1 & 1 \end{pmatrix}$ to each, and reattach them in order to get the nodal chart in Figure 1 c). Applying the shear matrix $\begin{pmatrix} 1 & 2 \\ 0 & 1 \end{pmatrix}$ to the nodal chart diagram in Figure 1 c), we get back the original nodal chart diagram Figure 1 a). This gives an integral affine isomorphism $\tau : B_0 \to B_1$.

**The Hamiltonian diffeomorphism.** Let $\pi : S^2 \times S^2 \times [0,1] \to B \times [0,1]$ be a lift of the nodal tangle $B \times [0,1]$ supported on $\mathscr{F}^{-1}((M-\varepsilon, M])$ ([3, Definition 2.30]). Since $\tau$ is an integral affine isomorphism identifying $B_0$ and $B_1$, [3, Theorem 2.26] gives a symplectomorphism $\psi \in \mathrm{Symp}(S^2 \times S^2)$ such that

$$\begin{array}{ccc}
S^2 \times S^2 & \xrightarrow{\psi} & S^2 \times S^2 \\
\downarrow{\pi_0} & & \downarrow{\pi_1} \\
(B_0, \mathfrak{N}_0) & \xrightarrow{\tau} & (B_1, \mathfrak{N}_1)
\end{array}$$

commutes outside of the set $\mathscr{F}^{-1}((M-\varepsilon, M])$. Let $x \in \mathscr{F}^{-1}([0, M-\varepsilon))$ in $B_0$, and $h = \mathscr{F}(x)$ its height. Then $\psi(\pi_0^{-1}(x)) = \pi_1^{-1}(\tau(x)) = \pi_0(2(M-h) \cdot x)$, where $t \cdot$ denotes the $\mathbb{R}$-action on $B_0 \setminus \mathscr{F}^{-1}(M)$ translating a point $x$ by integral affine distance $t$ in the clockwise direction along the level set $\mathscr{F}^{-1}(h)$. The level-set $\mathscr{F}^{-1}(h)$ has integral affine length $2w + 8(M-h)$, so for $x \in \mathscr{F}^{-1}([0, M-\varepsilon))$ with $\frac{w}{M-h}$ irrational we have $\psi^n(\pi_0^{-1}(x)) \cap \pi_0^{-1}(x) = \emptyset$ for all $n \geq 0$. The symplectic mapping class group of non-monotone $S^2 \times S^2$ is trivial, so $\psi \in \mathrm{Symp}_0(S^2 \times S^2) = \mathrm{Ham}(S^2 \times S^2)$.

## 2 The general construction

**The base space.** Let $(B_0, \mathfrak{N}_0)$ be a nodal integral affine surface as in [3, Proposition 3.17] with *canonical type* [3, Definition 3.23] $(H, M, \alpha_1, ..., \alpha_n)$ obtained by modifying a Delzant polygon by a nodal tangle. Recall that $(B_0, \mathfrak{N}_0)$ admits a natural nodal chart as in [3, Remark 3.18], illustrated in Figure 2 a): Let $\mathscr{F} : B_0 \to [0, M)$ be the height function on $B_0$ inherited from the Delzant polygon ([3, Section 3.1]). Let $\varepsilon > 0$. We may choose the nodes in the hat $H$ to have at least height $M - \varepsilon$. The level sets $\mathscr{F}^{-1}(h)$ with $h \in [0, M - \varepsilon)$ are straight lines and circles in $B_0$. We have $M = \max \mathscr{F}$, and we call $\Delta_M := \mathscr{F}^{-1}(M)$ the **ridge**

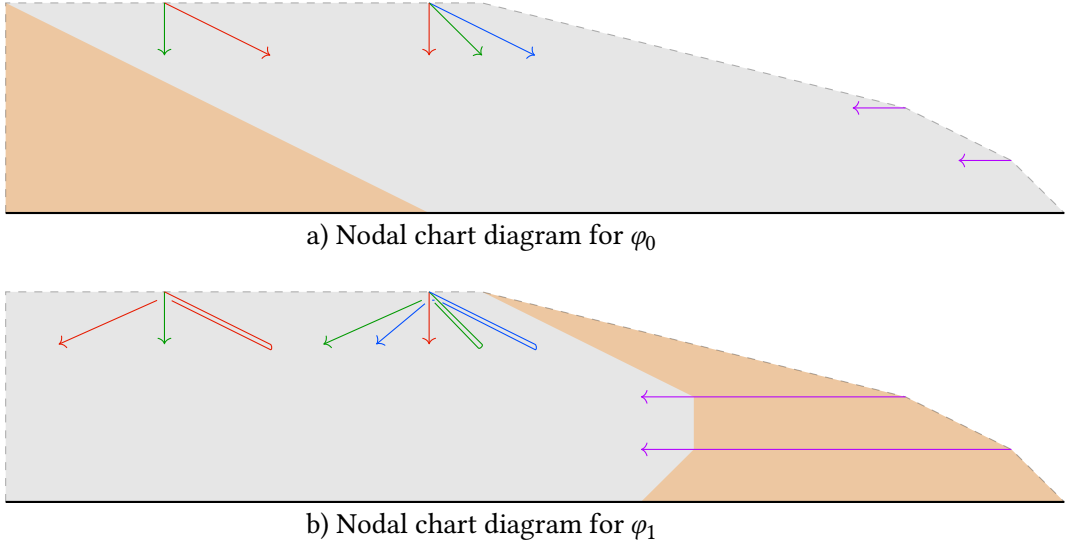

a) Nodal chart diagram for $\varphi_0$

b) Nodal chart diagram for $\varphi_1$

Figure 2: Two nodal charts for the same canonical type.

of $B_0$. The integral affine length of $\Delta_M$ is the *width* $w$ of the hat $H$ ([3, Definition 3.20]). The integral affine length of $\mathscr{F}^{-1}(h)$ is given by the function

$$g : \mathbb{R} \to \mathbb{R}$$
$$h \mapsto 2w - k(M - h) + \sum_{i=1}^{n} \min\{h - \alpha_i, 0\},$$

where $k$ is a constant depending on the hat class of $H$ as follows:

| Ⓑ | Ⓒ | Ⓓ | Ⓔ |
|---|---|---|---|
| $k = 8$ | $k = 8$ | $k = 7$ | $k = 6$ |

where we used the same notation as in Figure 3. Take a straight line segment $\ell$ connecting $\partial B_0$ to the ridge $\Delta_M$. Take $x \in B_0 \setminus \Delta_M$ and set $h = \mathscr{F}(x)$. Let $f : B_0 \setminus (\ell \cup \Delta_M) \to \mathbb{R}^2$ be the length of the straight line segment connecting $\mathscr{F}^{-1}(h) \cap \ell$ with $x$ along the level set $\mathscr{F}^{-1}(h)$ in the counter-clockwise direction.[1] Then $\varphi_0 : B_0 \setminus (\ell \cup \Delta_M) \to \mathbb{R}^2$ with $\varphi_0(x) = (f(x), \mathscr{F}(x))$ defines a nodal chart, as in Figure 2 a).

**The nodal tangle.** If $n > 0$, $B_0$ has at least one parked node. Let $\mathfrak{a}_n$ be a parked node with smallest height $\alpha_n = \mathscr{F}(\mathfrak{a}_n)$. Sliding $\mathfrak{a}_n$ once around the circle $\mathscr{F}^{-1}(\alpha_n)$ gives us the construction in [2, Remark 3.2] leading to a counterexample to the Lagrangian Poincaré recurrence conjecture.

Here we give a different approach that provides slightly different examples. Suppose that the hat $H$ has width $w > 0$. Then the hat $H$ contains the line segment $\Delta_M$. As can be seen from Figure 3, there are two types of ends of $\Delta_M$: either there are three nodes pointing to the end or two. Depending on the types of ends of $\Delta_M$, modify $H$ by a nodal tangle as in Figure 4. Let $\mathfrak{a}_1, ..., \mathfrak{a}_n$ be the parked nodes of $B_0$ with height $\mathscr{F}(\mathfrak{a}_i) = \alpha_i$. Further modify $B_0$ by a nodal tangle by sliding each parked node $\mathfrak{a}_i$ to the left by $2(M - \alpha_i)$. Denote by $(B_1, \mathfrak{N}_1)$ the resulting nodal integral affine surface, and by $(B \times [0, 1], \mathfrak{N})$ the nodal tangle connecting $(B_0, \mathfrak{N}_0)$ to $(B_1, \mathfrak{N}_1)$. We also get the natural nodal chart $\varphi_1 : (B_1, \mathfrak{N}_1) \to \mathbb{R}^2$ as in Figure 2 b). There is a integral affine isomorphism $\tau : (B_0, \mathfrak{N}_0) \to (B_1, \mathfrak{N}_1)$: The nodal chart diagrams in Figure 2 a) and b) of $(B_0, \mathfrak{N}_0)$ and $(B_1, \mathfrak{N}_1)$ are related by applying the shear matrix $\begin{pmatrix} 1 & 2 \\ 0 & 1 \end{pmatrix}$ to Figure 2 a) and cutting the orange triangle off Figure 2 a) and regluing it as indicated in Figure 2 b). We quickly verify that this indeed gives an integral affine isomorphism: For points in the interior of the nodal chart diagrams in Figure 2 this is obvious. For points in a neighbourhood of $\Delta_M$ this can be verified in a nodal chart as in Figure 4.

---

[1] Here we choose a orientation of $B_0$.

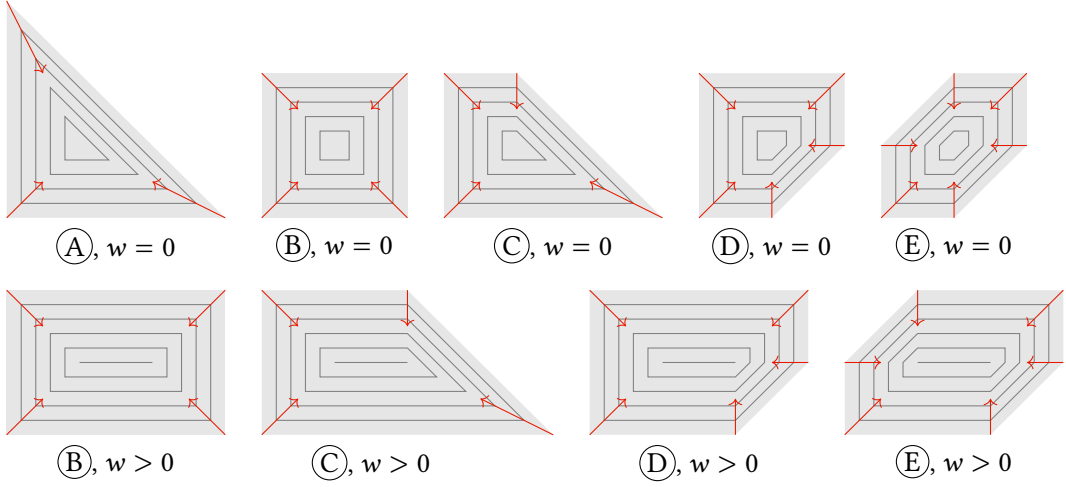

Figure 3: All possible $\varepsilon$-hats up to nodal tangle.

**The Hamiltonian diffeomorphism.** Let $U$ be an $\varepsilon$-neighbourhood in $[0, M]$ of the points $\{M, \alpha_1, ..., \alpha_n\}$. Then we have $\mathscr{F}(\pi_B(\mathfrak{N})) \subset U$. Let $\pi : X \times [0, 1] \to (B \times [0, 1], \mathfrak{N})$ be a lift of the nodal tangle $(B \times [0, 1], \mathfrak{N})$ supported on $\mathscr{F}^{-1}(U)$ ([3, Definition 2.30]). The symplectic manifold $X$ is determined by $(H, M, \alpha_1, ..., \alpha_n)$, see [3, Theorem 3.25]. Since $\tau$ is an integral affine isomorphism identifying $B_0$ and $B_1$, [3, Theorem 2.26] gives a symplectomorphism $\psi \in \text{Symp}(X)$ such that

$$\begin{array}{ccc} X & \xrightarrow{\psi} & X \\ \downarrow{\pi_0} & & \downarrow{\pi_1} \\ (B_0, \mathfrak{N}_0) & \xrightarrow{\tau} & (B_1, \mathfrak{N}_1) \end{array}$$

commutes outside of the set $\mathscr{F}^{-1}(U)$. Let $x \in \mathscr{F}^{-1}([0, M] \setminus U)$ in $B_0$, and $h = \mathscr{F}(x)$ its height. Then $\psi(\pi_0^{-1}(x)) = \pi_1^{-1}(\tau(x)) = \pi_0(2(M - h) \cdot x)$, where $t \cdot$ denotes the $\mathbb{R}$-action on $B_0 \setminus \Delta_M$ translating a point $x$ by integral affine distance $t$ in the clockwise direction along the level set $\mathscr{F}^{-1}(h)$. The level-set $\mathscr{F}^{-1}(h)$ has integral affine length $g(h)$, so for $x \in \mathscr{F}^{-1}([0, M] \setminus U)$ with $\frac{g(h)}{2(M-h)}$ irrational we have $\psi^n(\pi_0^{-1}(x)) \cap \pi_0^{-1}(x) = \emptyset$ for all $n \geq 0$. [1, Corollary 4.6] says that the symplectic mapping class group of $X$ is finite, thus $\psi^m \in \mathrm{Symp}_0(X) = \mathrm{Ham}(X)$ for some $m \geq 0$.

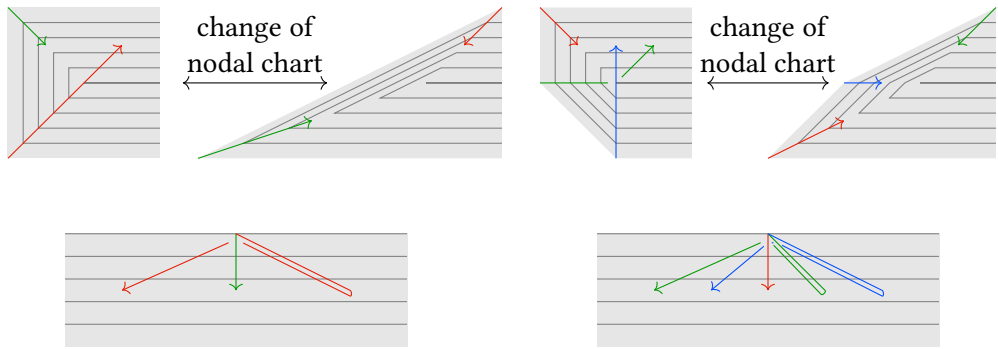

Figure 4: Nodal tangle transforming part of an $\varepsilon$-hat of positive width. On the left the nodal tangle for an end of $\Delta_M$ with two incoming nodes, on the right the one for an end with three incoming nodes. The top row shows the nodal tangle in a nodal chart as in Figure 3, the bottom row shows the same nodal tangle in a nodal chart as in Figure 2.